\documentclass{amsart}
\usepackage{amsmath,amssymb}
\usepackage{yhmath}

\newcommand{\bdis}{\begin{displaymath}}
\newcommand{\edis}{\end{displaymath}}
\newcommand{\be}{\begin{equation}}
\newcommand{\ee}{\end{equation}}
\newcommand{\mbb}{\mathbb}
\newcommand{\mcal}{\mathcal}

\newcommand{\vp}{\varphi} 
\newcommand{\vth}{\vartheta}

\newcommand{\zf}{\zeta\left(\frac{1}{2}+it\right)}

\newcommand{\FR}{\frac{x^n+y^n}{z^n}}

\DeclareMathOperator{\im}{Im}

\theoremstyle{definition}

\theoremstyle{remark}
\newtheorem{remark}[]{Remark}

\newtheorem*{mydef11}{{\bf Theorem 1}}

\newtheorem*{mydef12}{{\bf Theorem 2}}

\newtheorem*{mydef13}{{\bf Theorem 3}}

\newtheorem*{mydef41}{{\bf Corollary 1}}

\newtheorem*{mydef42}{{\bf Corollary 2}}

\newtheorem*{mydef43}{{\bf Corollary 3}}

\newtheorem*{mydef81}{{\bf Property 1}}

\newtheorem*{mydef82}{{\bf Property 2}}

\newtheorem*{mydef91}{{\bf Formula 1}}

\newtheorem*{mydef92}{{\bf Formula 2}}

\numberwithin{equation}{section}

\begin{document}

\title[Jacob's ladders and the equivalent of the Fermat-Wiles theorem \dots]{Jacob's ladders and the equivalent of the Fermat-Wiles theorem generated by the Hardy-Littlewood formula (1921)}

\author{Jan Moser}

\address{Department of Mathematical Analysis and Numerical Mathematics, Comenius University, Mlynska Dolina M105, 842 48 Bratislava, SLOVAKIA}

\email{jan.mozer@fmph.uniba.sk}

\keywords{Riemann zeta-function}

\begin{abstract}
In this paper we prove that the 105 years old Lemma 18 of Hardy and Littlewood from theirs fundamental paper \cite{2} is able to generate, for example, a continuum set of new $\zeta$-equivalents of the Fermat-Wiles theorem.  
\end{abstract}
\maketitle

\section{Introduction} 

\subsection{} 

Let us remind that Hardy and Littlewood have proved, see \cite{2}, Theorem B, the fundamental inequality 
\be \label{1.1} 
N_0(T+U)-N_0(T)>KU,\ T>T_0, 
\ee  
where by $N_0(T)$ denotes the number of zeros of the function $\zf$, $0<t<T$, and 
\be \label{1.2} 
U=T^a,\ a>\frac 12,\ K=K(a),\ T_0=T_0(a). 
\ee 
The proof of the estimate (\ref{1.1}) is based on the Lemma 18, see \cite{2}, p. 305: 
\be \label{1.3} 
J=\int_T^{T+U}\mcal{J}^2{\rm d}t=\pi\sqrt{2\pi}HU+\mcal{O}\left(\frac{U}{\ln T}\right), 
\ee  
where 
\be \label{1.4} 
\mcal{J}=\mcal{J}(t,H)=\int_t^{t+H}X(u){\rm d}u, 
\ee   
\be \label{1.5} 
X(t)=-\frac{t^{\frac 14}e^{\frac 14\pi t}}{\frac 14+t^2}\Xi(t), 
\ee 
and 
\be \label{1.6} 
T^a<U<T^b,\ a>\frac 12,\ 0<H<T^c. 
\ee   
Compare the first inequality in (\ref{1.6}) with the equality in (\ref{1.2}), where $c$ is positive and sufficiently small. 

Finally, Hardy and Littlewood have remarked\footnote{See \cite{2}, p. 315.} the following: \emph{As was observed in (5.5), we do not use the full force of Lemma 18. The complete Lemma, however, seems of considerable interest in itself, and may prove to be of service in the future.} 

\begin{remark}
In this paper we shall prove, after 105 years, this prophecy of Hardy and Littlewood that their Lemma 18 is able to generate, for example, the continuum set of new $\zeta$-equivalents of the Fermat-Wiles theorem. 
\end{remark} 

\subsection{} 

Next, let us remind the formula\footnote{See \cite{26}, p. 79.} 
\be \label{1.7} 
Z(t)=-2\pi^{\frac{1}{4}}\frac{\Xi(t)}{(\frac 14+t^2)|\Gamma(\frac 14+i\frac t2)|}. 
\ee 
Elimination of $\Xi(t)$ from (\ref{1.5}) and (\ref{1.7}) gives that 
\be \label{1.8} 
X(t)=\frac 12\pi^{-\frac 14}t^{\frac 14}e^{\frac 14\pi t}\left|\Gamma\left(\frac 14+i\frac t2\right)\right|Z(t), 
\ee 
and the Stirling formula 
\be \label{1.9} 
\left|\Gamma\left(\frac 14+i\frac t2\right)\right|\sim \sqrt{2\pi}e^{-\frac 14\pi t}\left(\frac t2\right)^{-\frac 14},\ t\to\infty 
\ee 
implies this result. 

\begin{remark}
It is the connection between the Hardy-Littlewood $X(t)$-function and the Riemann's $Z(t)$-function: 
\be \label{1.10} 
X(t)\sim \left(\frac{\pi}{2}\right)^{\frac 14}Z(t),\ t\to\infty. 
\ee 
These two functions shall appear simultaneously in our formulas. Furthermore, we use $\sqrt{2\pi^3}$ instead of $\pi\sqrt{2\pi}$ in the complicated formulas. 
\end{remark} 

\begin{remark}
The formula (\ref{1.5}) follows from the equality of the third and the fourth member in the last row on the page no. 289 in \cite{2}. However, in the third member, we must make up the minus sign since 
\bdis 
s=\frac 12+it \ \Rightarrow \ s(s-1)=-\left(\frac 14+t^2\right). 
\edis 
\end{remark} 

\subsection{} 

In this paper we obtain the following functional generated by the Hardy-Littlewood formula (\ref{1.3}): 
\be \label{1.11} 
\begin{split}
& \lim_{\rho\to\infty}\frac{1}{\rho}\left\{
\int_{[x\rho]^1}^{[x\rho+\frac{(2\pi^3)^{\frac{1}{2\alpha}}}{a}]^1}Z^2(t){\rm d}t
\right\}^{-\alpha}\times \\ 
& \int_{(x\rho)^{1/a}}^{(x\rho)^{1/a}+x\rho}\left[\int_t^{t+a^{-\alpha}(\ln x\rho)^\alpha}X(u){\rm d}u\right]^2{\rm d}t=x
\end{split}
\ee 
for every fixed $x>0,\ a>\frac 12,\ \alpha>0$, where 
\be \label{1.12} 
[G]^1=\vp_1^{-1}(G). 
\ee  
In the special case 
\be \label{1.13} 
x\to\FR,\ x,y,z,n\in\mbb{N},\ n\geq 3 
\ee 
of the Fermat's rationals we have the following result: the $\zeta$-condition 
\be \label{1.14} 
\begin{split}
	& \lim_{\rho\to\infty}\frac{1}{\rho}\left\{
	\int_{[\FR\rho]^1}^{[\FR\rho+\frac{(2\pi^3)^{\frac{1}{2\alpha}}}{a}]^1}Z^2(t){\rm d}t
	\right\}^{-\alpha}\times \\ 
	& \int_{(\FR\rho)^{1/a}}^{(\FR\rho)^{1/a}+\FR\rho}\left[\int_t^{t+a^{-\alpha}(\ln \FR\rho)^\alpha}X(u){\rm d}u\right]^2{\rm d}t\not=1 
\end{split}
\ee 
on the set of all Fermat's rationals expresses the next equivalent of the Fermat-Wiles theorem for every fixed $a>\frac 12$ and $\alpha>0$. 

\begin{remark}
Simultaneously, we have shown that there is a usage of above-mentioned Lemma 18 in a direction that is completely different as was the original usage of that statement. 
\end{remark} 

\subsection{} 

The fourth part of this work contains some remarks on A. A. Karatsuba usage of my results about the equation $Z'(t)=0$, see \cite{12}, in his paper \cite{4} dedicated to the 90th. anniversary of I. M. Vinogradov. Within these remarks we give also new \emph{extremely distanced}\footnote{Terminus used by I. M. Vinogradov.} limits of possibilities of the Weyl's  method of trigonometric sums in the corresponding problems. My remarks are based on I. M. Vinogradov analysis in the Introduction of his monograph \cite{27}.

\section{Jacob's ladders: notions and basic geometrical properties}  

\subsection{}

In this paper we use the following notions of our works \cite{18} -- \cite{22}: 
\begin{itemize}
\item[{\tt (a)}] Jacob's ladder $\vp_1(T)$, 
\item[{\tt (b)}] direct iterations of Jacob's ladders 
\bdis 
\begin{split}
	& \vp_1^0(t)=t,\ \vp_1^1(t)=\vp_1(t),\ \vp_1^2(t)=\vp_1(\vp_1(t)),\dots , \\ 
	& \vp_1^k(t)=\vp_1(\vp_1^{k-1}(t))
\end{split}
\edis 
for every fixed natural number $k$, 
\item[{\tt (c)}] reverse iterations of Jacob's ladders 
\be \label{2.1}  
\begin{split}
	& \vp_1^{-1}(T)=\overset{1}{T},\ \vp_1^{-2}(T)=\vp_1^{-1}(\overset{1}{T})=\overset{2}{T},\dots, \\ 
	& \vp_1^{-r}(T)=\vp_1^{-1}(\overset{r-1}{T})=\overset{r}{T},\ r=1,\dots,k, 
\end{split} 
\ee   
where, for example, 
\be \label{2.2} 
\vp_1(\overset{r}{T})=\overset{r-1}{T}
\ee  
for every fixed $k\in\mbb{N}$ and every sufficiently big $T>0$. We also use the properties of the reverse iterations listed below.  
\be \label{2.3}
\overset{r}{T}-\overset{r-1}{T}\sim(1-c)\pi(\overset{r}{T});\ \pi(\overset{r}{T})\sim\frac{\overset{r}{T}}{\ln \overset{r}{T}},\ r=1,\dots,k,\ T\to\infty,  
\ee 
\be \label{2.4} 
\overset{0}{T}=T<\overset{1}{T}(T)<\overset{2}{T}(T)<\dots<\overset{k}{T}(T), 
\ee 
and 
\be \label{2.5} 
T\sim \overset{1}{T}\sim \overset{2}{T}\sim \dots\sim \overset{k}{T},\ T\to\infty.   
\ee  
\end{itemize} 

\begin{remark}
	The asymptotic behaviour of the points 
	\bdis 
	\{T,\overset{1}{T},\dots,\overset{k}{T}\}
	\edis  
	is as follows: at $T\to\infty$ these points recede unboundedly each from other and all together are receding to infinity. Hence, the set of these points behaves at $T\to\infty$ as one-dimensional Friedmann-Hubble expanding Universe. 
\end{remark}  

\subsection{} 

Let us remind that we have proved\footnote{See \cite{22}, (3.4).} the existence of almost linear increments 
\be \label{2.6} 
\begin{split}
& \int_{\overset{r-1}{T}}^{\overset{r}{T}}\left|\zf\right|^2{\rm d}t\sim (1-c)\overset{r-1}{T}, \\ 
& r=1,\dots,k,\ T\to\infty,\ \overset{r}{T}=\overset{r}{T}(T)=\vp_1^{-r}(T)
\end{split} 
\ee 
for the Hardy-Littlewood integral (1918), \cite{1}: 
\be \label{2.7} 
J(T)=\int_0^T\left|\zf\right|^2{\rm d}t. 
\ee  

For completeness, we give here some basic geometrical properties related to Jacob's ladders. These are generated by the sequence 
\be \label{2.8} 
T\to \left\{\overset{r}{T}(T)\right\}_{r=1}^k
\ee 
of reverse iterations of the Jacob's ladders for every sufficiently big $T>0$ and every fixed $k\in\mbb{N}$. 

\begin{mydef81}
The sequence (\ref{2.8}) defines a partition of the segment $[T,\overset{k}{T}]$ as follows 
\be \label{2.9} 
|[T,\overset{k}{T}]|=\sum_{r=1}^k|[\overset{r-1}{T},\overset{r}{T}]|
\ee 
on the asymptotically equidistant parts 
\be \label{2.10} 
\begin{split}
& \overset{r}{T}-\overset{r-1}{T}\sim \overset{r+1}{T}-\overset{r}{T}, \\ 
& r=1,\dots,k-1,\ T\to\infty. 
\end{split}
\ee 
\end{mydef81} 

\begin{mydef82}
Simultaneously with the Property 1, the sequence (\ref{2.8}) defines the partition of the integral 
\be \label{2.11} 
\int_T^{\overset{k}{T}}\left|\zf\right|^2{\rm d}t
\ee 
into the parts 
\be \label{2.12} 
\int_T^{\overset{k}{T}}\left|\zf\right|^2{\rm d}t=\sum_{r=1}^k\int_{\overset{r-1}{T}}^{\overset{r}{T}}\left|\zf\right|^2{\rm d}t, 
\ee 
that are asymptotically equal 
\be \label{2.13} 
\int_{\overset{r-1}{T}}^{\overset{r}{T}}\left|\zf\right|^2{\rm d}t\sim \int_{\overset{r}{T}}^{\overset{r+1}{T}}\left|\zf\right|^2{\rm d}t,\ T\to\infty. 
\ee 
\end{mydef82} 

It is clear, that (\ref{2.10}) follows from (\ref{2.3}) and (\ref{2.5}) since 
\be \label{2.14} 
\overset{r}{T}-\overset{r-1}{T}\sim (1-c)\frac{\overset{r}{T}}{\ln \overset{r}{T}}\sim (1-c)\frac{T}{\ln T},\ r=1,\dots,k, 
\ee  
while our eq. (\ref{2.13}) follows from (\ref{2.6}) and (\ref{2.5}).  

\section{New $\zeta$-functional and corresponding $\zeta$-equivalent of the Fermat-Wiles theorem generated by the Hardy-Littlewood Lemma 18} 

\subsection{} 

In the case 
\be \label{3.1} 
U=T^a,\ H=\ln^\alpha T,\ a>\frac 12,\ \alpha>0 
\ee 
which fulfils the assumptions of the Hardy-Littlewood Lemma 18\footnote{See (\ref{1.6}), comp. (\ref{1.2}).}, where, of course, 
\bdis 
H=\ln^\alpha T<T^c,\ T\to\infty
\edis 
for every fixed $x>0$ and $c>0$, we have the following formula\footnote{See (\ref{1.3}) and (\ref{1.4}).} 
\be \label{3.2} 
\begin{split}
& \int_T^{T+T^a}\left[\int_t^{t+\ln^\alpha T}X(u){\rm d}u\right]{\rm d}t=\\ 
& \{1+o(1)\}\sqrt{2\pi^3}T^a\ln^\alpha T,\ T\to\infty. 
\end{split}
\ee 
Now, we use the substitution 
\be \label{3.3} 
T=\tau^{\frac 1a},\ \{T\to +\infty\} \ \Leftrightarrow \ \{\tau\to +\infty\}
\ee 
in (\ref{3.2}) that gives this result 
\be \label{3.4} 
\begin{split}
& \int_{\tau^{\frac 1a}}^{\tau^{\frac 1a} +\tau}\left[\int_t^{t+a^{-\alpha}\ln^\alpha\tau}X(u){\rm d}u\right]{\rm d}t= \\ 
& \{1+o(1)\}\sqrt{2\pi^3}a^{-\alpha}\tau\ln^\alpha\tau,\ \tau\to\infty. 
\end{split}
\ee 
In the right-hand side of (\ref{3.4}) we next use our formula\footnote{See \cite{23}, (5.1), $k=1$.} in the form 
\be \label{3.5} 
\ln\tau=\frac{1+o(1)}{2l}\int_{[\tau]^1}^{[\tau+2l]^1}Z^2(t){\rm d}t,\ l>0 
\ee 
that gives us the following formula 
\be \label{3.6} 
\begin{split}
& \int_{\tau^{\frac 1a}}^{\tau^{\frac 1a} +\tau}\left[\int_t^{t+a^{-\alpha}\ln^\alpha\tau}X(u){\rm d}u\right]{\rm d}t= \\ 
& \{1+o(1)\}\frac{\sqrt{2\pi^3}}{(2la)^\alpha}\tau\left(\int_{[\tau]^1}^{[\tau+2l]^1}Z^2(t){\rm d}t\right)^2,\ \tau\to\infty. 
\end{split}
\ee 
Further, we choose for $2l$ the value 
\be \label{3.7} 
2l_1=\frac{(2\pi^3)^{\frac{1}{2\alpha}}}{a}, 
\ee  
and after this we use the substitution 
\be \label{3.8} 
\tau=x\rho, \ \{\tau\to+\infty\} \ \Leftrightarrow \ \{\rho\to+\infty\}, x>0 
\ee 
in the formula (\ref{3.6}), $l=l_1$, that gives us the following functional. 

\begin{mydef11}
\be \label{3.9} 
\begin{split}
& \lim_{\rho\to\infty}\frac{1}{\rho}\left\{
\int_{[x\rho]^1}^{[x\rho+\frac{(2\pi^3)^{\frac{1}{2\alpha}}}{a}]^1}Z^2(t){\rm d}t
\right\}^{-\alpha}\times \\ 
& \int_{(x\rho)^{1/a}}^{(x\rho)^{1/a}+x\rho}\left[\int_t^{t+a^{-\alpha}(\ln x\rho)^\alpha}X(u){\rm d}u\right]^2{\rm d}t=x
\end{split}
\ee 
for every fixed $x>0$, $a>\frac 12$ and $\alpha>0$. 
\end{mydef11} 

\subsection{} 

In the special case of the Fermat's rationals, see (\ref{1.13}), we obtain the following. 

\begin{mydef41}
\be \label{3.10} 
\begin{split} 
	& \lim_{\rho\to\infty}\frac{1}{\rho}\left\{
\int_{[\FR\rho]^1}^{[\FR\rho+\frac{(2\pi^3)^{\frac{1}{2\alpha}}}{a}]^1}Z^2(t){\rm d}t
\right\}^{-\alpha}\times \\ 
& \int_{(\FR\rho)^{1/a}}^{(\FR\rho)^{1/a}+\FR\rho}\left[\int_t^{t+a^{-\alpha}(\ln \FR\rho)^\alpha}X(u){\rm d}u\right]^2{\rm d}t=\FR 
\end{split} 
\ee 
for every fixed Fermat's rational and every fixed $a>\frac 12$ and $\alpha>0$. 
\end{mydef41} 

Consequently, the following Theorem holds true. 

\begin{mydef12}
The $\zeta$-condition 
\be \label{3.11} 
\begin{split} 
	& \lim_{\rho\to\infty}\frac{1}{\rho}\left\{
	\int_{[\FR\rho]^1}^{[\FR\rho+\frac{(2\pi^3)^{\frac{1}{2\alpha}}}{a}]^1}Z^2(t){\rm d}t
	\right\}^{-\alpha}\times \\ 
	& \int_{(\FR\rho)^{1/a}}^{(\FR\rho)^{1/a}+\FR\rho}\left[\int_t^{t+a^{-\alpha}(\ln \FR\rho)^\alpha}X(u){\rm d}u\right]^2{\rm d}t\not= 1 
\end{split} 
\ee 
on the set of all Fermat's rationals expresses the new $\zeta$-equivalent of the Fermat-Wiles theorem for every fixed $a>\frac 12$ and every fixed $\alpha>0$. 
\end{mydef12} 

\begin{remark}
Since the symbol 
\bdis 
\{(a,\alpha):\ a>\frac 12, \alpha>0\}
\edis 
ranges over the continuum set, then the formula (\ref{3.11}) represents continuum set of new $\zeta$-equivalents of the Fermat-Wiles theorem. 
\end{remark} 

\section{Remarks on A. A. Karatsuba's use of my results on eq. $Z'(t)=0$ in his paper dedicated to the 90th. anniversary of  I. M. Vinogradov} 

Let us remind that we have obtained new new type of functionals and corresponding $\zeta$-equivalents of the Fermat-Wiles theorem in our paper \cite{24}. These results are based on our asymptotic formula (1981) for the function $Z'(t)$, see \cite{12}, (2), (16) - the second formula. 

\begin{remark}
This fact alone justifies few remarks about a manner in which A. A. Karatsuba has used my original results about the equation $Z'(t)=0$ in his paper \cite{4}. 
\end{remark}
 
The structure of our further treatise is as follows. First, we remind: 
\begin{itemize}
	\item[(a)] Our basic formulas concerning the equation $Z'(t)=0$, \cite{12}. 
	\item[(b)] Our original results following the above-mentioned formulas.  
	\item[(c)] The scepticism of I. M. Vinogradov about the possibilities of the method of trigonometric sums, see \cite{27}, p. 13. 
	\item[(d)] Some remarks on A. A. Karatsuba's use of our results in his paper \cite{4}. 
\end{itemize} 

\subsection{} 

Let us remind that Riemann has defined the real-valued function 
\be \label{4.1.1}  \tag{4.1.1}
Z(t)=e^{i\vth(t)}\zf, 
\ee 
where\footnote{See \cite{25}, (35), (44), (62), comp. \cite{26}, p. 98.} 
\be \label{4.1.2} \tag{4.1.2} 
\begin{split}
& \vth(t)=-\frac t2\ln\pi+\im\ln\Gamma\left(\frac 12+i\frac t2\right)= \\ 
& \frac t2\ln\frac{t}{2\pi}-\frac t2-\frac{\pi}{8}+\mcal{O}\left(\frac 1t\right).  
\end{split}
\ee  
Next, let 
\be \label{4.1.3} \tag{4.1.3} 
S(a,b)=\sum_{0<a\leq n<b\leq 2a}n^{it},\ b\leq\sqrt{\frac{t}{2\pi}}
\ee 
be an elementary trigonometric sum. We have obtained in our paper \cite{12}, (2) and (16), new asymptotic formulae: if 
\be \label{4.1.4} \tag{4.1.4} 
|S(a,b)|< A(\Delta)\sqrt{a}t^\Delta,\ 0<\Delta\leq \frac 16, 
\ee  
then 
\be \label{4.1.5} \tag{4.1.5} 
\begin{split}
& \sum_{T\leq\bar{t}_{2\nu}\leq T+H}Z'(\bar{t}_{2\nu})=-\frac{1}{4\pi}H\ln^2\frac{T}{2\pi}+\mcal{O}(T^\Delta\ln^2T), \\ 
& \sum_{T\leq\bar{t}_{2\nu+1}\leq T+H}Z'(\bar{t}_{2\nu+1})=\frac{1}{4\pi}H\ln^2\frac{T}{2\pi}+\mcal{O}(T^\Delta\ln^2T),
\end{split}
\ee 
where the sequence 
\be \label{4.1.6} \tag{4.1.6} 
\{\bar{t}_\nu\}_{\nu=1}^\infty 
\ee 
is defined by the condition 
\be \label{4.1.7} \tag{4.1.7} 
\vth(\bar{t}_\nu)=\pi\nu+\frac{\pi}{2}, 
\ee  
and further 
\be \label{4.1.8} \tag{4.1.8} 
H=T^{\Delta}\psi(T), 
\ee  
where the function $\psi(T)$ is a function increasing (arbitrarily slowly) to $+\infty$. 

\begin{remark}
Here, it is sufficient to put $\psi(T)=\ln T$, i.e. 
\be \label{4.1.9} \tag{4.1.9} 
H=T^\Delta\ln T. 
\ee 
\end{remark} 

The proof of our formulae (\ref{4.1.5}) is then based on the following\footnote{See \cite{12}, (10), and the proof (17) -- (39).}

\begin{mydef91}
\be \label{4.1.10} \tag{4.1.10} 
Z'(t)=-2\sum_{n<\sqrt{\frac{t}{2\pi}}}\frac{1}{\sqrt{n}}(\vth-\ln n)\sin(\vth-t\ln n)+\mcal{O}(t^{-1/4}\ln t), 
\ee 
\end{mydef91} 
that we have transformed into 
\begin{mydef92}
\be \label{4.1.11} \tag{4.1.11} 
\begin{split}
& Z'(t)=-2\sum_{n<P_0}\frac{1}{\sqrt{n}}\ln\frac{P_0}{n}\sin(\vth-t\ln n)+\mcal{O}(T^{-1/4}\ln T), \\ 
& t\in [T,T+H],\ H\in (0,\sqrt[4]{T}],\ P_0=\sqrt{\frac{T}{2\pi}}. 
\end{split}
\ee 
\end{mydef92} 

\begin{remark}
The factor 
\be \label{4.1.12} \tag{4.1.12} 
\ln\frac{P_0}{n},\ n<P_0 
\ee 
in the amplitudes of the Riemann's oscillators is of the prime importance in the formula (\ref{4.1.11}) since this one enables us to remove small denominators\footnote{Comp. \cite{12}, (57) -- (64).} in the course of estimating the corresponding sums. 
\end{remark} 

\subsection{} 

Here we enumerate our results. 

\begin{itemize}
\item[(A)] The main {\bf Theorem} of our paper \cite{12} follows from the formulae (\ref{4.1.5}): If 
\be \label{4.2.1} \tag{4.2.1} 
|S(a,b)|<A(\Delta)\sqrt{a}t^\Delta,\ 0<\Delta\leq \frac 16, 
\ee 	
then there is a root of odd order of the equation 
\be \label{4.2.2} \tag{4.2.2} 
Z'(t)=0 
\ee 
within the interval\footnote{Comp. (\ref{4.1.8}) and (\ref{4.1.9}).} 
\be \label{4.2.3} \tag{4.2.3} 
(T,T+T^{\Delta}\ln T), 
\ee  
see \cite{12}, (2) -- (4).  
\item[(B)] For example, in the case of Kolesnik's exponent, see \cite{6} 
\be \label{4.2.4} \tag{4.2.4} 
\Delta=\frac{35}{216}+\epsilon
\ee 
we have the following. 
\begin{mydef41}
The interval 
\be \label{4.2.5} \tag{4.2.5} 
(T,T+T^{35/216+\epsilon})
\ee 
contains the root of the odd order of the equation (\ref{4.2.2}), where $\epsilon>0$ is arbitrarily small, \cite {12}, (6). 
\end{mydef41} 
\begin{remark}
A. Kolesnik's exponent $\frac{35}{216}$ represents only $2.8\%$-improvement of the exponent $\frac{1}{6}$. 
\end{remark}
\item[(C)] On the Lndel\" of hypothesis we have\footnote{See \cite{3}, p. 89.} 
\be \label{4.2.6} \tag{4.2.6} 
|S(a,b)|<A(\epsilon)\sqrt{a}t^\epsilon, 
\ee 
for every sufficiently small $\epsilon>0$, i.e. we have the following. 
\begin{mydef42}
On the Lindel\" of hypothesis the interval 
\be \label{4.2.7} \tag{4.2.7} 
(T,T+T^\epsilon\ln T)
\ee 
contains the root of odd order of the equation (\ref{4.2.2}). 
\end{mydef42} 
\begin{remark}
This means that the Lindel\" of hypothesis gives a $100\%$-improvemenst of the Kolesnik's exponent in (\ref{4.2.5}) as well as it gives a $100\%$-improvement to all possible further improvements of the Kolesnik's exponent. 
\end{remark} 
\item[(D)] Next, let $N'_0(T_0)$ stands for the number of the roots of the equation (\ref{4.2.2}) in the interval $(0,T)$. We have obtained the following. 
\begin{mydef43}
On the Lindel\" of hypothesis we have the following estimate 
\be \label{4.2.8} \tag{4.2.8} 
N_0'(T+T^\tau)-N_0'(T)>A(\tau,\epsilon)T^{\tau-\epsilon},\ 0<\epsilon<\tau, 
\ee  
where $\tau$ is an arbitrarily small fixed number.  
\end{mydef43} 
\end{itemize} 

\subsection{} 

Next, we remind the I. M. Vinogradov's scepticism on possibilities of the method of trigonometric sums (H. Weyl's sums). 

Namely, I. M. Vinogradov has analysed in the Introduction of his monograph \cite{27} the possibilities of the method of trigonometric sums in the problem of estimation the remainder term $R(N)$ in the asymptotic formula\footnote{See \cite{27}, p. 13.}
\be \label{4.3.1} \tag{4.3.1} 
\pi(N)-\int_2^{N}\frac{{\rm d}x}{\ln x}=R(N), 
\ee  
where we have: 
\begin{itemize}
	\item[(a)] the original Vall\' e - Poussin's result 
	\be \label{4.3.2} \tag{4.3.2} 
	R(N)=\mcal{O}(Ne^{-c\sqrt{\ln N}}), 
	\ee  
	\item[(b)] H. Weyl's result 
	\be \label{4.3.3} \tag{4.3.3} 
	R(N)=\mcal{O}(Ne^{-c_1\sqrt{\ln N\ln\ln N}}), 
	\ee 
	\item[(c)] I. M. Vinogradov's result 
	\be \label{4.3.4} \tag{4.3.4} 
	\begin{split}
	& R(N)=\mcal{O}(Ne^{-c_3\lambda(N)\ln^{0.6}N}), \\ 
	& \lambda(N)=(\ln\ln N)^{-0.2}. 
	\end{split}
	\ee 
\end{itemize}

Under the Riemann hypothesis, however, the next estimate follows (H. von Koch (1901), \cite{5}) 
\be \label{4.3.5} \tag{4.3.5} 
R(N)=\mcal{O}(\sqrt{N}\ln N), 
\ee  
and this one is extremely distanced from (\ref{4.3.4}). 

In connection with this fact, I. M. Vinogradov made this remark: 

\emph{Obviously, it is very hard to make any essential progress in solution of the problem to find the order of $R(N)$ term willing to find $R(N)=\mcal{O}(N^{1-\gamma}), \gamma=0.000001$ by making use of only some improvements of the H. Weyl's estimates and without making use of further important progress in the theory of zeta-function.} 

In the following subsections we give some remarks about a manner of using (or not using) of our original results about the equation $Z'(t)=0$ in the A. A. Karatsuba's paper \cite{4}. 

\subsection{} 

It is noticed by A. A. Karatsuba, see \cite{4}, p. 305: \emph{Here I prove a new theorem on the zeros of $Z^{(k)}(t)$ that, in the cases $k=0,1$, improves corresponding J. Moser's result.} 

\begin{remark}
A. A. Karatsuba's asserion is, however, incorrect since: 
\begin{itemize}
	\item[(a)] Our Theorem (the main result), that is completely ignored in the paper \cite{4}, is general for the parameters $\Delta\in (0,1/6]$ and its structure is as follows: if (\ref{4.2.1}) then (\ref{4.2.3}) with the corresponding property. 
	\item[(b)] A. A. Karatsuba has selected from our results in \cite{12} only an example of our general theorem i.e. only Corollary 1, see (\ref{4.2.4}) and (\ref{4.2.5}) but, without the assumption of this corollary. 
	\item[(c)] Next, A. A. Karatsuba did not prove the estimate 
	\be \label{4.4.1} \tag{4.4.1} 
	|S(a,b)|<A_1\sqrt{a}t^{1/12}, 
	\ee 
	and therefore, his result about the property of the interval 
	\be \label{4.4.2} \tag{4.4.2} 
	(T,T+T^{1/12}\ln T) 
	\ee  
	is not directly connected with ours (Moser and Kolesnik) Corollary 1. 
	\item[(d)] Of course, our Corollary 2\footnote{See (\ref{4.2.7}) and Remark 11.} gives immediately $100\%$-improvement (on the Lindel\" of hypothesis) of A. A. Karatsuba's exponent $\frac{1}{12}$ in (\ref{4.4.2}) and what is more important: $100\%$-improvement if all such $\Delta$ that fulfils the condition (\ref{4.2.1}). 
\end{itemize}
\end{remark} 

\subsection{} 

Our Theorem, together with the Corollary 2, has also an other goal. 

\begin{remark}
Namely, this one gives some test on the possibilities of improvements of estimates of the elementary trigonometric sum in our problem. This test represents, of course, an analogue of I. M. Vinogradov analysis, see subsection 4.3. 
\end{remark}

Example: 
\begin{itemize}
	\item[(a)] Let us assume that the following estimate 
	\be \label{4.5.1} \tag{4.5.1} 
	|S(a,b)|<B_1\sqrt{a}t^{0.01}
	\ee 
	holds true, probably, it is not achieved yet. 
	\item[(b)] This assumption immediately generates, for example, the following question: when it will be achieved the estimate 
	\be \label{4.5.2} \tag{4.5.2} 
	\begin{split}
	& |S(a,b)|<B_2\sqrt{a}t^\omega, \\ 
	& \omega = 0.000\ 000 \ 001?
	\end{split}
	\ee 
\end{itemize} 

\begin{remark}
Wee see, that the limit for the improvements of the elementary trigonometric sum that gives the Lindel\" of hypothesis is \emph{extremely distanced} from all already-reached estimates. This is the full analogy of I. M. Vinogradov case, comp. (\ref{4.3.4}) and (\ref{4.3.5}). 
\end{remark}
 
\subsection{} 

Now, we examine the influence of the Riemann hypothesis on the limit of improvements of the interval 
\be \label{4.6.1} \tag{4.6.1} 
(T,T+T^\Delta\ln T) 
\ee 
with the corresponding property. Namely, J. E. Littlewood has obtained, see \cite{7}, p. 237, under the Riemann hypothesis the estimate 
\be \label{4.6.2} \tag{4.6.2} 
\gamma''-\gamma'<\frac{A}{\ln\ln\gamma'},\ \gamma'\to\infty
\ee 
for every consecutive zeros 
\be \label{4.6.3} \tag{4.6.3} 
\frac 12+i\gamma',\ \frac 12+i\gamma'',\ \gamma'<\gamma'' 
\ee  
of the function $\zeta(s)$. And from this one gets immediately (in our context) that every interval 
\be \label{4.6.4} \tag{4.6.4} 
\left(\gamma',\gamma'+\frac{A}{\ln\ln\gamma'}\right) 
\ee 
contains an odd-order zero of the function $Z'(t)$. 

\begin{remark}
On Riemann hypothesis the interval 
\be \label{4.6.5} \tag{4.6.5} 
\left(T, T+\frac{2A}{\ln\ln T}\right),\ T\to\infty
\ee 
contains the zero of odd-order of the function $Z'(t)=0$. Of course, to reach the limit (\ref{4.6.5}) by means of improvements of the intervals (\ref{4.6.1}) on the base of elementary trigonometric sum $S(a,b)$ only is probably an image from Lewis Caroll's Alice's Adventures in Wonderland \dots 
\end{remark}

\subsection{} 

Next, it follows from the Remark 15 the following. 

\begin{remark}
On Riemann hypothesis we have the following properties: 
\begin{itemize}
	\item[(a)] The interval 
	\be \label{4.7.1} \tag{4.7.1} 
	\left( T, T+\frac{1}{\ln\ln\ln T}\right)
	\ee 
	contains at least 
	\be \label{4.7.2} \tag{4.7.2} 
	\left[ B_3\frac{\ln\ln T}{\ln\ln\ln T}\right]=k_1(T) 
	\ee 
	zeros of odd order of the function 
	\be \label{4.7.3} \tag{4.7.3} 
	Z^{(1)}(t),\ t\in \left( T, T+\frac{1}{\ln\ln\ln T}\right). 
	\ee  
	\item[(b)] And, of course, the interval (\ref{4.7.1}) contains a zero of odd order of the function 
	\be \label{4.7.4} \tag{4.7.4} 
	Z^{(k_1(T))}(t) 
	\ee 
	by the elementary Rolle's theorem only. 
\end{itemize} 

Next, let us remind A. A. Karatsuba's result, see \cite{4}, p. 63, Corollary 2: for every $\epsilon>0$ there is $k_0(\epsilon)>0$ such that the interval 
\be \label{4.7.5} \tag{4.7.5} 
(T, T+T^\epsilon)
\ee 
contains odd-order zero of the function 
\be \label{4.7.6} \tag{4.7.6} 
Z^{(k)}(t),\ t\in (T,T+T^\epsilon),\ k>k_0(\epsilon). 
\ee 
\end{remark} 

\begin{remark}
If we compare (\ref{4.7.5}) with (\ref{4.7.1}) and also (\ref{4.7.6}) and (\ref{4.7.4}) then we obtain the following mappings 
\be \label{4.7.7} \tag{4.7.7} 
\begin{split}
& T^\epsilon \searrow \frac{1}{\ln\ln\ln T}, \\ 
& k\nearrow \left[ B_3\frac{\ln\ln T}{\ln\ln\ln T}\right]=k_1(T). 
\end{split}
\ee 
We see that the Riemann hypothesis generates the non-reachable limits 
\be \label{4.7.8} \tag{4.7.8} 
\frac{1}{\ln\ln\ln T},\ \left[ B_3\frac{\ln\ln T}{\ln\ln\ln T}\right]=k_1(T) 
\ee 
for the improvements of A. A. Karatsuba's (\ref{4.7.5}) and (\ref{4.7.6}) based on the improvements of the estimates of Weyl's sums. This conclusion is compatible also with I. M. Vinogradov analysis in subsection 4.3. 
\end{remark} 

Consequently, we give the following. 

\begin{mydef13}
The Riemann hypothesis implies the next \emph{limit estimate}: within every infinitesimal interval 
\be \label{4.7.9} \tag{4.7.9} 
\left( T, T+\frac{1}{\ln\ln\ln T}\right),\ T\to\infty 
\ee  
are packed zeros of odd order of all the functions 
\be \label{4.7.10} \tag{4.7.10} 
\begin{split}
& Z^{(k)}(t),\ t\in \left( T, T+\frac{1}{\ln\ln\ln T}\right), \\ 
& k=1,2,\dots, \left[ B_3\frac{\ln\ln T}{\ln\ln\ln T}\right], 
\end{split}
\ee 
see (\ref{4.7.1}) -- (\ref{4.7.4}). 
\end{mydef13} 

\subsection{} 

A. A. Karatsuba has written in his paper \cite{4}, p. 52, these words in the connection with his Lemma 5: \emph{I did not obtain in the literature similar assertion and quote it \dots} 
\begin{remark} 
This statement does not correspond with the reality since:  
\begin{itemize}
	\item[(a)] His Lemma 5 gives, in the case $k=0$ (the second equation), the following 
	\be \label{4.8.1} \tag{4.8.1} 
	Z^{(1)}(t)=-2\sum_{n\leq \sqrt{\frac{t}{2\pi}}}\frac{\vth'(t)-\ln n}{\sqrt{n}}\sin(\vth(t)-t\ln n)+\mcal{O}(t^{-1/4}\ln^2t), 
	\ee   
	and this one is almost identical with my formulae (\ref{4.1.10}) that gives only a bit better error term 
	\bdis 
	\mcal{O}(t^{-1/4}\ln t). 
	\edis 
	\item[(b)] Our Formula 2, see (\ref{4.1.11}), contains new and basic factor 
	\be \label{4.8.2} \tag{4.8.2} 
	\ln\frac{P_0}{n}=\left(\ln\frac{P_0}{n}\right)^1. 
	\ee 
	\item[(c)] I have sent the paper \cite{12} to Acta Arithmetica in the beginning of 1979 (30. 01. 1979) and simultaneously to A. A. Karatsuba.  
	\item[(d)] Of course, our paper \cite{12} contains also the complete proof of the formula (\ref{4.1.10}), see \cite{12}, (17) -- (39). 
\end{itemize}
\end{remark}  

\begin{remark}
That is: 
\begin{itemize}
	\item[(a)] A. A. Karatsuba's Lemma 5 represents the generalization of my original formula (\ref{4.1.10}). 
	\item[(b)] His corresponding formula in \cite{4}, p. 59, represents only generalization\footnote{Comp. (\ref{4.8.2}).} 
	\bdis 
	\left(\ln\frac{P_0}{n}\right)^1\to \dots \to \left(\ln\frac{P_0}{n}\right)^k 
	\edis 
	of our original formula (\ref{4.1.11}). 
\end{itemize}
\end{remark} 

\subsection{} 

On the basis of essential improvements of the classical results of Hardy and Littlewood (1918) and (1921) and also of the proof of Titchmarsh hypothesis (1934) I was invited to take part on the international conference in Moscow organized in honour of the 90th birthday of I. M. Vinogradov in the period Sept. 14 -- 19, 1981. 

\begin{remark}
This conference has joined my research with that of A. A. Karatsuba, namely: 
\begin{itemize}
	\item[(a)] I have presented a talk that has listed main results obtained in the period 1976 -- 1981, see \cite{17}, comp. \cite{8} -- \cite{16} in this paper. 
	\item[(b)] A. A. Karatsuba published , see \cite{4}, the first work from his series of papers in the direction of my research, I have sent him the separates as well as the published manuscripts of the papers accepted for publication. 
\end{itemize}
\end{remark} 

\begin{remark}
Our scientific positions were as follows: 
\begin{itemize}
	\item[(c)] A. A. Karatsuba was full Professor and Member of the Steklov Institute in Moscow. 
	\item[(d)] I was a lector without scientific degree (1966 -- 1990), studying individually the works of Hardy and Littlewood, Titchmarsh, Ingham, and, of course, Selberg on the Riemann zeta-function.  
\end{itemize}
\end{remark} 

\begin{remark}
Only after the fall of communism in Czechoslovakia in 1989 I have obtained directly the highest academic degree D.Sc. on the basis of my dissertation \emph{The evolution of the Titchmarsh discreet method}, pp. I -- XIV, 1 -- 236 (in Russian), unpublished. I am grateful for this to the kind initiative of Prof. Stefan Schwarz, D. Sc., co-founder of the theory of semigroups and the Member of the editorial board of the journal \emph{Semigroups forum}. 
\end{remark}

I would like to thank Michal Demetrian for his moral support of my study of Jacob's ladders.


\begin{thebibliography}{29}
\bibitem{1} 
G.H. Hardy, J.E. Littlewood, Contribution to the theory of the Riemann zeta-function and the theory of the distribution of Primes, Acta Math. 41 (1), 119 -- 196, (1918).  
\bibitem{2} 
G.H. Hardy, J.E. Littlewood, The zeros of Riemann's zeta-function on the critical line, Math. Z. 10, 283--317, (1921). 
\bibitem{3} 
A. A. Karatsuba, \emph{Basic analytic number theory}, Springer Verlag, 1992. 
\bibitem{4} 
A. A. Karatsuba, On the distance between adjacent zeros of the Riemann zeta-function on the critical line, Trudy Math. Inst. Steklov, 157, 49--63, (1981), (in Russian). 
\bibitem{5} 
H. von Koch, Sur la distribution des nombres premiers, Acta Math. 24, 159--181, (1901). 
\bibitem{6} 
A. Kolesnik, On the order of the Dirichlet $L$-function, Pacific J. Math., 82, 159--182, (1979). 
\bibitem{7} 
J. E. Littlewood, Two notes on the Riemann zeta-function, Proc. Camb. Phil. Soc., 22, 234-242, (1924). 
\bibitem{8} 
J. Moser, On certain sum in the theory of the Riemann zeta-function, Acta Arith. 31, 31-43, (1976), (in Russian). 
\bibitem{9} 
J. Moser, On a theorem of Hardy and Littlewood in the theory of the Riemann zeta-function, Acta Arith. 31, 45-51, (1976), (in Russian). 
\bibitem{10} 
J. Moser, Supplement to \cite{9}, Acta Arith. 35, 403--404, (1979), (in Russian). 
\bibitem{11} 
J. Moser, Proof of the Titchmarsh conjecture in the theory of the Riemann zeta-function, Acta Arith. 36, 147--156,  (1980), (in Russian). 
\bibitem{12} 
J. Moser, On the roots of the equation $Z'(t)=0$, Acta Arith. 40, 79--89, (1981), (in Russian). 
\bibitem{13} 
J. Moser, Correction to \cite{8}--\cite{10}, Acta Arith. 40, 97--107, (1981), (in Russian). 
\bibitem{14} 
J. Moser, On a lemma of Hardy and Littlewood i the theory of the Riemann zeta-function, Acta Math. Univ. Comeniane, 42(43), 7--26, (1983), (in Russian). 
\bibitem{15} 
J. Moser, Sharpening of the Hardy-Littlewood theorem on the density of zeros of the function $\zf$, Acta Math. Univ. Comeniane, 42(43), 41-50, (1983), (in Russian). 
\bibitem{16} 
J. Moser, Improvement of the theorem of Hardy and Littlewood on density of the zeros of the function $\zf$, Acta Arith. 43, 21--47, (1983), (in Russian). 
\bibitem{17} 
J. Moser, Some consequences of the Riemann-Siegel formula, Trudy Math. Inst. Steklov, 163, 183--185, (1983), (in Russian). 
\bibitem{18}
J. Moser,
`Jacob's ladders and almost exact asymptotic representation of the Hardy-Littlewood integral`,
Math. Notes 88, (2010), 414-422, arXiv: 0901.3937. 
\bibitem{19}
J. Moser,
`Jacob's ladders, the structure of the Hardy-Littlewood integral and some new class of nonlinear integral equations`,
Proc. Steklov Inst. 276 (2011), 208-221, arXiv: 1103.0359. 
\bibitem{20}
J. Moser, Jacob's ladders, reverse iterations and new infinite set of $L_2$-orthogonal systems generated by the Riemann $\zf$-function, arXiv: 1402.2098v1.  
\bibitem{21} 
J. Moser, Jacob's ladders, existence of almost linear increments of the Hardy-Littlewood integral and new types of multiplicative laws, arXiv: 2304.09267. 
\bibitem{22} 
J. Moser, Jacob's ladders, almost linear increments of the Hardy-Littlewood integral (1918) and their relations to the Selberg formula (1946) and the Fermat-Wiles theorem, arXiv: 2312.12085. 
\bibitem{23} 
J. Moser, Jacob's ladders, $\zeta$-transformation of the Fourier orthogonal system (2014) and new infinite set of $\zeta$-equivalents of the Fermat-Wiles theorem, arXiv: 2507.06724v1. 
\bibitem{24} 
J. Moser, Jacob's ladders, our asymptotic formulae (1981) connected with the equation $Z'(t)=0$ as a base for new $\zeta$-equivalents of the Fermat-Wiles theorem and some other results, arXiv: 2509.17626. 
\bibitem{25} 
C. L. Siegel, \" Uber Nachlass zur analytischen Zahlentheorie, Quellen und Studien zur Geschichte der Mat. Astr. und Physik, Abt. B: Studien, 2, 45--80, (1932).  
\bibitem{26} 
E. C. Titchmarsh, \emph{The theory of the Riemann zeta-function}, Oxford, At the Clarendon Press, 1951. 
\bibitem{27} 
I. M. Vinogradov, \emph{The method of trigonometric sums in the number theory}, NAUKA, Moscow, 1981.   
\end{thebibliography}
\end{document}